\documentclass[preprint,10pt]{article}

\usepackage{amssymb}
\usepackage{color}
 \usepackage{amsthm}
\usepackage{amsmath}

\bibliography{alpha.bst}

\newtheorem{thm}{Theorem}[section]

\newtheorem{lem}[thm]{Lemma}
\newtheorem{cor}[thm]{Corollary}

\newtheorem{rem}[thm]{Remark}
\newcommand{\eproof}{\rule{0,2cm}{0,2cm}}

\begin{document}

\begin{center}

{\huge {\bf On fractional Duhamel's}}

\smallskip

{\huge {\bf principle and its applications}}

\vspace{1cm}

{\large {\bf Sabir Umarov}}

\vspace{.5cm}

{Tufts University, Medford, MA 02155, USA\\}

\end{center}

\begin{abstract}
The classical Duhamel principle, established nearly 200 years ago by
Jean-Marie-Constant Duhamel, reduces the Cauchy problem for an
inhomogeneous partial differential equation to the Cauchy problem
for the corresponding homogeneous equation. Duhamel's principle is
not applicable in the case of fractional order differential
equations. In this paper we formulate and prove fractional
generalizations of this famous principle directly applicable to a
wide class of fractional order differential-operator equations.

\end{abstract}

\label{}
\section{Introduction}

Let $X$ be a reflexive Banach space and $A: {\mathcal{D}}
\rightarrow X$ a closed linear operator with a domain ${\mathcal{D}}
\subset X$. In Section 2 we will introduce a Frech\'et type
topological vector space $\mathrm{Exp}_{A,G}(\mathrm{X})$ (and its
dual $\mathrm{Exp}^{'}_{A^{\ast}, G^{\ast}}(\mathrm{X^{\ast}})$),
where $G$ is an open subset of the complex plain $\mathbb{C}.$
This space represents a modification of the space of entire
functions with finite exponential type
\cite{Djrbashian93,Du82,Levin,Radino} and its abstract versions. We
also introduce a functional calculus in the form $f(A),$ where $f$
is an analytic function defined on $G.$ The function $f$ is called
the symbol of the operator $f(A).$

The goal of this  paper is to generalize Duhamel's principle for the
Cauchy problem for general inhomogeneous fractional distributed
order differential-operator equations of the form
\begin{align} \label{eq1}
L^{\Lambda}[u] &\equiv \int_0^{\mu}
f(\alpha,A)D_{\ast}^{\alpha}u(t)d \Lambda (\alpha)
= h(t),\quad t>0,\\
\label{cauchy1}
u^{(k)}(0)&= \varphi_{k}, \quad
k = 0, ..., m-1,
\end{align}
where  $\mu \in (m-1,m];$ $h(t)$ and $\varphi_{k}, \,\,k = 0,...,
m-1,$ are given $X$-valued
vector-functions; $f(\alpha,A)$ is a family of operators with the
symbol $f(\alpha,z)$ continuous in the variable $\alpha \in
[0,\mu],$ and analytic in the variable $z \in G \subset
{\mathbb{C}};$ $\Lambda$ is a finite measure defined on $[0,\mu];$
and $D_{\ast}^{{\alpha}}$ is the operator of fractional
differentiation of order $\alpha$ in the sense of Caputo-Djrbashian
(see, for example, \cite{Ca67,GM97}), i. e.
\begin{equation*}
 D_{\ast}^\alpha g(t)=
      \begin{cases}
      \frac{1}{\Gamma(n - \alpha)} \int_{0}^{t}
        \frac{g^{(n)}(\tau)d\tau}{(t-\tau)^{\alpha-n+1}}, & \text{if  $n-1
        < \alpha < n, \, n \in \mathbb{N}$,}\\
       g^{(n)}(t) \equiv \frac{d^n}{dt^n}g(t), & \text{if $\alpha=n \in \{0\} \cup \mathbb{N}$.}
       \end{cases}
\end{equation*}
Hereafter the integrals  are understood in the sense of Bochner if
$g(t)$ is a vector-function with values in some topological-vector
space for each fixed $t.$

The classical Duhamel principle is not applicable in the case of
fractional order differential equations. Its modification combined
with some integral transformations can reduce the Cauchy problem for
inhomogeneous equation to the Cauchy problem for homogeneous
equation. However, this two step process becomes cumbersome for
complex equations containing many terms with fractional operators.
In this paper we formulate and prove {\em fractional generalizations
of Duhamel's principle}  applicable directly to the Cauchy problem
for inhomogeneous fractional order differential-operator equations,
which reduce them to the  Cauchy problem for corresponding
homogeneous equations. In the particular case of fractional order
partial differential equations with a single "fractional" term in
the equation (\ref{eq1}), a fractional analog of Duhamel's principle
was obtained in \cite{US06,US07}.

Fractional order differential equations are useful and appropriate
mathematical apparatus for modeling problems with memory, and
interest in this subject has grown substantially during the last a
few decades. For instance, probability density functions of a wide
class of non-Gaussian diffusion processes satisfy fractional order
governing equations with space and time fractional order
differential operators (see \cite{AM03,MMM,HKU,MK00} and references
therein). Inhomogeneous fractional order differential equations
appear naturally describing the influence of an external force or
memory effects. In the study of diffusion processes in complex
heterogeneous media with several distinct diffusion modes, even
without an external force, the function $h(t)$ embodies memory of
the past \cite{LH,UmarovSteinberg}.

There is extensive literature on the Cauchy problem for integer
order abstract differential-operator equations (see, e.g.
\cite{DalKrein,VasPis}). The first order evolution equations
$u^{'}(t)=Au(t)$ in the spaces of abstract exponential
vector-functions of a finite type, $Exp_A(X)$ (and in more general
bornological spaces) were studied in \cite{Radino}.
In the case of integer $\alpha_{k}, \,\,k = 1,..., m$, the Cauchy
problem for pseudo-differential and differential-operator equations
with analytic symbols or with symbols having singularities was
studied, for example, in \cite{Du81,Um86,TV89}, and multi-point
value problems in \cite{Pt84,SU99,Um97,Um98,Um92}. What concerns
fractional order differential-operator equations, Kochubei
\cite{Koc89} studied existence and uniqueness of a solution to the
abstract Cauchy problem $D_{\ast}^{\alpha}u(t)=Au(t), \, u(0)=u_0,$
with Caputo-Djrbashian fractional derivative for $0<\alpha<1$ and a
closed operator $A$ with a dense domain ${\cal{D}}(A)$ in a Banach
space. El-Sayed \cite{E-S95}, Bazhlekova \cite{Ba2001} investigated
the Cauchy problem for $0<\alpha<2.$ In the more general case of
$\alpha
>0,$ Gorenflo et. al. \cite{GLZ99} studied existence of solutions in
Roumieu-Beurling and Gevrey classes. Kostin \cite{Kos93} proved
correctness of the abstract initial value problem (Cauchy type
problem) $D_+^{\alpha}u(t)=Au(t), \, D_+^{\alpha-k}u(0)=\varphi_k,
\, k=1,...,m,$ for $\alpha \in (m-1,m),$ and with the
Rieman-Liouville derivative $D_+^{\alpha}.$  For more information
about recent results on the Cauchy problem for abstract fractional
differential-operator equations, we refer the reader to
\cite{Ba2001,EidKoch2004,KilbasST}; and for a recent mathematical
treatment of the distributed fractional order differential equations
to papers \cite{Kochubei2008,MSh,UG2005}.

The paper is organized as follows. In Section 2 we recall the
classic Duhamel principle, and the basic spaces of elements used in
this paper. Since we formulate a fractional Duhamel principle in the
abstract case, we introduce a topological-vector space on which the
corresponding operators act. In Section 3 we formulate the main
result, namely an abstract fractional analog of Duhamel's principle
and discuss some of its applications.

\section{Preliminaries}
\subsection{Fractional order derivatives.}

For a function $g$ defined on $[0,\infty)$,
under some integrability conditions the {\it fractional integral of
order $\beta$} with  \emph{terminal points} $\tau$ and $t,$ is
defined as
\[
_{\tau}J^{\beta}g(t)=\frac{1}{\Gamma(\beta)}\int_\tau^{t}
(t-s)^{\beta -1}g(s)ds,
\]
where $\Gamma(\cdot)$ is Euler's gamma-function. Obviously, if
$\beta =n$ then $_{\tau}J^n$ is the $n$-fold integral of $f$ over
the interval $[\tau,t].$ By convention, $_{\tau}J^0 f(t)=f(t),$ i.e.
$_{\tau}J^{0}$ coincides with the identity operator. In the notation
we do not indicate the upper terminal point $t$, since in the
current paper it is always $t.$

Further, let $m$ be a positive integer number. We denote by $_\tau
D_{+}^{\alpha}, \, m-1 < \alpha < 1,$ the fractional derivative of
order $\alpha $ in the sense of Riemann-Liouville, which is defined
as
$$_\tau D_{+}^{\alpha}g(t) = \frac{1}{\Gamma(m-\alpha)}\frac{d^m}{dt^m}
\int_{\tau}^{t}\frac{g(s)ds}{(t-s)^{\alpha+1-m}} , \,\,m-1 < \alpha
< m,$$ and $_\tau D_{+}^{0}g(t)=g(t), \,\, _\tau D_{+}^{m}g(t)=
g^{(m)}(t)$. Between this fractional derivative and the
\emph{Caputo-Djrbashian} derivative there is the following
relationship \cite{GM97}:
\begin{equation}
\label{relation1} _\tau D_{+}^{\alpha}g(t) = \,
\,_{\tau}D_{*}^{\alpha}g(t) + \sum_{k=0}^{m-1}\frac{g^{(k)}(\tau)}
{\Gamma(k-\alpha+1)}(t-\tau)^{k-\alpha}, \, t>0.
\end{equation}
In the particular case of $0<\alpha<1$ one has
\begin{equation}
\label{relation2} _\tau D_{+}^{\alpha}g(t) =
\,_{\tau}D_{*}^{\alpha}g(t) +
g(\tau)\frac{(t-\tau)^{-\alpha}}{\Gamma(1-\alpha)}, \, t>0.
\end{equation}
If $g(\tau)=0,$ then one obtains the equality $_\tau
D_{+}^{\alpha}g(t) = \,_{\tau}D_{\ast}^{\alpha}g(t).$
 Alternative representations via the fractional integral
 are:
$$\,_{\tau}D_{+}^{\alpha}g(t)=\frac{d^m}{dt^m} \,_{\tau}J^{m-\alpha}g(t) ~~~ \mbox{and} ~~~
\,_{\tau}D_{\ast}^{\alpha}g(t)= \,_{\tau}J^{m-\alpha}\frac{d^m
g(t)}{dt^m}.$$
We omit the lower terminal point $\tau$ if $\tau=0,$
writing simply $D_{+}^{\alpha}, \, D_{*}^{\alpha}$ or $J^{\alpha}.$
Recall that for the Laplace transform of $D_{+}^{\alpha}g(t)$ and
$D_{\ast}^{\alpha}g(t),$ where $\alpha \in (m-1,m],$ the following
formulas are valid \cite{GM97}:
\begin{equation}
{\cal{L}}[D_{+}^{\alpha}g](s) = s^{\alpha} {\cal{L}}[g](s) -
\sum_{k=0}^{m-1} \frac{d^k}{dt^k} \left(J^{(m-\alpha)} g
\right)_{(t=0+)} s^{m-1-k},
\end{equation}
and
\begin{equation}
\label{laplace} {\cal{L}}[D_{\ast}^{\alpha}g](s) = s^{\alpha}
{\cal{L}}[g](s) - \sum_{k=0}^{m-1} g^{(k)} (0+) s^{\alpha -1 -k}.
\end{equation}
Here ${\cal{L}}[g](s)$ denotes the Laplace transform of $g.$

\subsection{An operator calculus}
\label{fcalculus}

In this section we recall some necessary facts about abstract spaces
of analytic elements of finite exponential type, and an operator
calculus defined on it. See for details \cite{Um97,Um98}.

Let ${X}$ be a reflexive Banach space with a norm $\|v\|, \, v \in
X$. Let $A$ be a closed linear operator with a domain ${\cal{D}}\,
(A)$ dense in ${X}$ and a spectrum $\sigma (A) \subset
\cal{C}.$ 
 Assume that $\sigma(A)$ is not empty and is not bounded.

We will develop an operator calculus $f(A)$ for analytic functions
$f(\lambda)$ in an open domain $G \subset \cal{C}.$ If the domain of
analyticity of $f,$ $G$ contains $\sigma(A)$ then
\begin{equation}\label{1}
f(A) = \int_{\nu} {\cal{R}} (\zeta, A) f(\zeta) d \zeta,
\end{equation}
where $\nu$ is a contour in $G$ containing $\sigma(A),$ and
$\cal{R}(\zeta, A), ~ \zeta \in \cal{C} \setminus \sigma(A),$ is the
resolvent operator of $A.$ However, if $f$ has singular points in
the spectrum of $A,$ then $f(A)$  can not be defined through the
integral \eqref{1}.

Assume that $G$ is any open set in $\cal{C}$ not necessarily
containing $\sigma(A).$ Further, let $0 < r \leq +\infty$ and $\nu <
r$. Denote by $\mathrm{Exp}_{A,\,\nu}(\mathrm{X})$ the set of
elements $v \in \cap_{k\ge 1}{\mathcal{D}}(A^k)$ satisfying the
inequalities $\|A^{k}v \| \leq C \nu^{k} \|v\|$ for all $k = 1, 2,
...,$ with a constant $C>0$ not depending on $k$. An element $v \in
\mathrm{Exp}_{A,\,\nu}({X})$ is said to be a vector of exponential
type $\nu$ \cite{Radino}. A sequence of elements $v_{n},\,\,n=1, 2,
...,$ is said to converge to an element $v_{0}$ in
$\mathrm{Exp}_{A,\,\nu}(\mathrm{X})$ iff:

 1) All the elements $v_{n}$ are vectors of exponential
 type $\nu < r,$ and

 2) $\|v_{n} - v_{0}\| \rightarrow 0, \,\, n \rightarrow
 \infty$.

\noindent Obviously, $\mathrm{Exp}_{A,\,\nu_1}(\mathrm{X}) \subset
\mathrm{Exp}_{A,\,\nu_2}(\mathrm{X}),$ if $\nu_1 < \nu_2.$
 Let \,\,$\mathrm{Exp}_{A,\,r}(\mathrm{X})$\,\, be the inductive limit of
 spaces \,\,$\mathrm{Exp}_{A,\,\nu}(\mathrm{X})$ when $\nu \rightarrow
 r$. For basic notions of topological vector spaces including inductive and
 projective limits we refer the reader to \cite{Robertson}. Set $A_{\lambda} = A-\lambda I,$ where $\lambda \in G,$
 and denote $\mathrm{Exp}_{A, r, \lambda}(\mathrm{X}) = \{ u_{\lambda} \in X:
u_{\lambda} \in \mathrm{Exp}_{A_{\lambda}, \, r}(\mathrm{X})\},$
with the induced topology. Finally, for
 arbitrary $G \subset \sigma(A),$ denote by
 $\mathrm{Exp}_{A,\,G}(\mathrm{X})$ the space whose elements
 are the locally finite sums of  elements in
 $\mathrm{Exp}_{A,\,r,\,\lambda}(\mathrm{X})$,\,\,$\lambda \in G, \,\,r < dist(\lambda, \partial G),$
 with the corresponding topology.
Namely, any $u\in \mathrm{Exp}_{A,\,G}(\mathrm{X})$ has a
representation $u=\sum_{\lambda}u_{\lambda}$ with a finite sum.
 It is clear, that $\mathrm{Exp}_{A,\,G}(\mathrm{X})$
 is a subspace of the space of vectors of exponential type if
 $r < +\infty,$ and coincides with it if $r = +\infty$.
 Moreover, $\mathrm{Exp}_{A,\,G}(\mathrm{X})$ is an abstract analog
 of the space $\Psi_{G,\,p}(R^1)$ introduced in \cite{Um97}, where
 $A = -i\frac{d}{dx}$, \,\,$G \subseteq R^{1}$,\,\,
 $\mathrm{X}$\,$ = L_p(R^1),\,\,1 < p < \infty$. In the case $A = -i\frac{d}{dx}$,\,\,
 $\mathrm{X}$\,$= L_2(R^1),$ the corresponding space was studied in
 \cite{Du81}.

 Further, let $f(\lambda)$ be an analytic function on $G.$ An
 arbitrary element $u \in \mathrm{Exp}_{A,\,G}(\mathrm{X})$ is
 represented as a finite sum
 $u = \sum_{\lambda \in G} u_{\lambda},\,\,\,u_{\lambda} \in
  \mathrm{Exp}_{A,\,R,\,\lambda}(\mathrm{X}).$
  Then for  $u \in \mathrm{Exp}_{A,\,G}(\mathrm{X})$ the
  operator $f(A)$ is defined by the formula
\begin{equation}
\label{calc1} f(A)u = \sum_{\lambda \in G}
f_{\lambda}(A)u_{\lambda},
\end{equation}
  where
\begin{equation}
\label{calc2} f_{\lambda}(A)u_{\lambda} = \sum_{n=0}^{\infty}
  \frac{f^{(n)}(\lambda)}{n!} A_{\lambda}^{n}u_{\lambda}.
  \end{equation}
In other words, each $f_{\lambda}$ represents locally $f$ in a
neighborhood of $\lambda \in G,$ and for $u_{\lambda}$ the operator
$f_{\lambda}(A)$ is well defined.

Additionally assume that there exists a one-parameter family of
bounded invertible operators $U_{\lambda}: {X}\rightarrow {X}$ such
that
\begin{equation}
\label{ulambda} AU_{\lambda} - U_{\lambda}A = \lambda U_{\lambda},
\,\,\lambda \in \sigma (A).
\end{equation}
For example, if $X=L_2 \equiv L_2(R)$ and $A = -i\frac{d}{dx}: L_2
\rightarrow L_2$ is the operator of differentiation with domain
${\mathcal{D}} (A) = \{ v \in L_2: A v \in L_2 \}$, then for the
operator $U_{\lambda}: v(x) \rightarrow e^{i\lambda x} v(x) $ we
have
\[
AU_{\lambda}v(x)=-\frac{d}{dx}(e^{i\lambda x}v(x)) = \lambda
e^{i\lambda x}v(x) - i e^{i\lambda x}\frac{dv}{dx} = \]
\[ \lambda U_{\lambda} v(x) + U_{\lambda} A v(x),
\]
obtaining (\ref{ulambda}). Condition (\ref{ulambda}) indicates a
shift of the spectrum of operator $A$ to $\lambda.$ This is seen
from the relationship $A-\lambda I = U_{\lambda}AU_{\lambda}^{-1},$
which follows from (\ref{ulambda}) multiplying by $U_{\lambda}^{-1}$
from the right. It follows from the latter that
\begin{equation}
\label{a^n} (A-\lambda I)^n = U_{\lambda} A^n U_{\lambda}^{-1},
\end{equation}
for all $n = 1,2,....$

Let $X^{\ast}$ denote the conjugate of $X,$ and $A^{\ast}: X^{\ast}
\to X^{\ast}$ be the  operator conjugate to $A$. Further, denote by
$\mathrm{Exp}^{'}_{A^{\ast}, G^{\ast}}(\mathrm{X^{\ast}})$ the space
of linear continuous functionals defined on $\mathrm{Exp}_{A,
G}(\mathrm{X}),$ with respect to weak convergence. Specifically, a
sequence $u_m^{\ast} \in \mathrm{Exp}^{'}_{A^{\ast},
G^{\ast}}(\mathrm{X^{\ast}})$ converges to an element $u^{\ast} \in
\mathrm{Exp}^{'}_{A^{\ast}, G^{\ast}}(\mathrm{X^{\ast}})$ if for all
$v \in \mathrm{Exp}_{A, G}(\mathrm{X})$ the convergence
$<u_m^{\ast}-u,v> \to 0$ holds as $m \to \infty.$  For an analytic
function $f^{w}$ defined on $G^{\ast}=\{z \in {\cal{C}}: \bar{z} \in
G\}$, we define {\it a weak extension} of $f(A)$ as follows:
\[
<f^{w}(A^{\ast})u^{\ast},v>=<u^{\ast}, f(A)v>, ~~ \forall v \in
\mathrm{Exp}_{A, G}(\mathrm{X}).
\]

\begin{lem}
Let $X$ be a reflexive Banach space and $A$ be a closed operator
defined on ${\cal{D}}(A) \subset X.$ Let $f$ be an analytic function
defined on a open connected set $G \subset \cal{C}.$ Then the
following mappings are well defined and continuous:
\begin{enumerate}
\item
$f(A): \mathrm{Exp}_{A, G}(\mathrm{X}) \rightarrow \mathrm{Exp}_{A,
G}(\mathrm{X}),$
\item
$f^{w}(A^{\ast}): \mathrm{Exp}^{'}_{A^{\ast},
G^{\ast}}(\mathrm{X^{\ast}}) \rightarrow \mathrm{Exp}^{'}_{A^{\ast},
G^{\ast}}(\mathrm{X^{\ast}}).$

\end{enumerate}

\end{lem}

{\em Proof.} Notice that $f(A)$ maps $\mathrm{Exp}_{A,
G}(\mathrm{X})$ into itself. Let $u \in \mathrm{Exp}_{A,
G}(\mathrm{X})$ has a representation $u=\sum_{\lambda}u_{\lambda}, ~
u_{\lambda} \in \mathrm{Exp}_{A_{\lambda}, \nu}(\mathrm{X}).$ Then
for $f(A)u$ defined in \eqref{calc1}, one has the following estimate

\begin{equation}
\label{estimate1} \|A^{k}_{\lambda} f_{\lambda}(A) u_{\lambda} \|
\le \sum_{n=0}^{\infty} \frac{|f^n(\lambda)|}{n!}\|(A-\lambda
I)^nA_{\lambda}^ku_{\lambda}\| \leq C \nu^{k} \|u_{\lambda}\|.
\end{equation}
 with some $\nu <r.$ It follows that $f_{\lambda}(A)u_{\lambda} \in \mathrm{Exp}_{A_{\lambda},
\nu}(\mathrm{X})$ with the same $\nu,$ and $f(A)u \in
\mathrm{Exp}_{A, G}(\mathrm{X}).$ The estimate (\ref{estimate1})
also implies continuity of the mapping $f(A)$ in the topology of
$\mathrm{Exp}_{A_{\lambda}, G}(\mathrm{X}).$

Now assume that a sequence $u_n^{\ast} \in
\mathrm{Exp}^{'}_{A^{\ast}, G^{\ast}}(\mathrm{X^{\ast}})$ converges
to $0$ in the weak topology of $\mathrm{Exp}^{'}_{A^{\ast},
G^{\ast}}(\mathrm{X^{\ast}}).$ Then for arbitrary $u \in
\mathrm{Exp}_{A, G}(\mathrm{X})$ we have
\[
<f^{w}(A^{\ast}) u_n^{\ast}, u> = <u_n^{\ast}, f(A)u>
=<u_n^{\ast},v,>
\]
where $v=f(A)u \in \mathrm{Exp}_{A, G}(\mathrm{X})$ due to the first
part of the proof. Hence, $f^{w}(A^{\ast})x_n^{\ast} \rightarrow 0,$
as $n \rightarrow \infty,$ in the weak topology of
$\mathrm{Exp}^{'}_{A^{\ast}, G^{\ast}}(\mathrm{X^{\ast}}).$ \eproof

\begin{rem}
 It is not hard to see that the above constructions are valid with
corresponding specifications in the case of operators with discrete
spectrum as well. Note that in this case the space $\mathrm{Exp}_{A,
G}(\mathrm{X})$ consists of the root lineals of eigenvectors
corresponding to the part of $\sigma(A)$ with nonempty intersection
with $G.$ If the spectrum $\sigma(A)$ is empty then an additional
exploration is required for solution spaces to be non-trivial (for
details see, {\cite{Du91}}).
\end{rem}

 As is shown in \cite{Um98}, the space $\mathrm{Exp}_{A,\,G}(\mathrm{X})$
 is invariant with respect to the action of an operator $f(A)$ and
 this operator acts continuously.
 \vspace {0.2 cm}

\subsection{Two lemmas}

The following two lemmas will be usful in proofs of theorems in
Section 3.
\begin{lem} \label{lemma1}
Let $h(t)$ be a continuous differentiable function. Then the
equation $J^{\alpha}u(t)=h(t), \, t>0,$ where $0<\alpha<1,$ has a
unique continuous solution given by the formula
\begin{equation} \label{Abel}
u(t)=D_+^{1-\alpha}h(t), \, t>0.
\end{equation}
\end{lem}

Lemma \ref{lemma1} is essentially the well-known result on a
solution of Abel's integral equation of first kind. See
\cite{GM97,SKM93} for the proof.

\begin{lem} \label{lemma2}
Suppose $v(t,\tau)$ is a vector-function in a Banach space $X,$
defined for all $t \ge \tau \ge 0$ and $k$ times differentiable with
respect to the variable $t.$ Let $u(t)=\int_0^t v(t,\tau)d\tau.$
Then
\begin{equation} \label{l2}
\frac{d^k}{dt^k}u(t)=\sum_{j=0}^{k-1}\frac{d^j}{dt^j}\big[\frac{\partial^{k-1-j}}{\partial
t^{k-1-j}}v(t,\tau)_{|_{\tau=t}}\big]+\int_0^t \frac{\partial^k
}{\partial t^k}v(t,\tau)d\tau.
\end{equation}
\end{lem}

\textit{Proof.} For a fixed $t>0$ and small $h$ one can easily
verify that
\begin{align}
\frac{u(t+h)-u(t)}{h}&=\frac{1}{h}\big(\int_0^{t+h}v(t+h,\tau)d\tau-\int_0^tv(t,\tau)d\tau\big)\notag\\
&=\frac{1}{h}\int_t^{t+h}v(t,\tau)d\tau
+\int_0^{t}\frac{v(t+h)-v(t)}{h}d\tau \notag\\
&+\int_t^{t+h}\frac{v(t+h)-v(t)}{h}d\tau. \label{es1}
\end{align}
Making use of the mean value theorem (in the integral form), we
obtain
\begin{align}
&\|\frac{1}{h}\int_t^{t+h}v(t,\tau)d\tau - v(t,t)\|\le
C_1\|v(t,\tau_{\ast})-v(t,t)\|, \, t<\tau_{\ast}<t+h, \label{es2}\\
&\|\int_0^{t}\frac{v(t+h)-v(t)}{h}d\tau - \int_0^t\frac{\partial
v(t,\tau)}{\partial t}d\tau\| \le C_2 |h|, \label{es3}\\
&\|+\int_t^{t+h}\frac{v(t+h)-v(t)}{h}d\tau\| \le C_3 |h|,
\label{es4}
\end{align}
where constants $C_1,\, C_2,$ and $C_3$ do not depend on $h.$ Now,
letting $h \to 0,$ estimates \eqref{es2}-\eqref{es4} and equation
\eqref{es1} imply the following formula:
\begin{equation} \label{l22}
\frac{d}{dt}u(t)=v(t,t)+\int_0^t \frac{\partial}{\partial
t}v(t,\tau)d\tau.
\end{equation}
Formula \eqref{l2} follows from \eqref{l22} by differentiation
repeatedly. \eproof

\subsection{Classical Duhamel's principle.} 
Duhamel's principle was formulated first for the Cauchy problem for
 second order linear inhomogeneous differential equations. Let \,$B=B(x,\frac{\partial}{\partial t},
 D_x),$\, where
$D_x=(\frac{\partial}{\partial x_1},...,\frac{\partial}{\partial
x_n}),$ be a linear differential operator with coefficients not
depending on \,$t$, and containing temporal derivatives of order not
higher than 1. Consider the Cauchy problem
\begin{equation}
\label{classic1} \frac{\partial ^{2}u}{\partial t^{2}}(t,x) + B
u(t,x) = f(t,x), \quad t > 0,\,\,x \in R^{n},
\end{equation}
with homogeneous initial conditions
\begin{equation}
\label{classic2} u(0,x) = 0,\quad \frac{\partial u}{\partial t}(0,x)
= 0.
\end{equation}
Let a sufficiently smooth function \,$v(t, \tau, x),\,\,t \geq
\tau,\,\,\tau \geq 0,\,\, x \in R^{n},$\, be for \,$t > \tau$\, a
solution of the homogeneous equation
$$
\frac{\partial^{2}v}{\partial t^{2}}(t, \tau, x) + B v(t, \tau, x) =
0,
$$
satisfying the following conditions:
$$
v(t, \tau, x)|_{t=\tau} = 0,\quad \frac{\partial v}{\partial t}(t,
\tau, x)|_{t=\tau}
 = h(\tau, x).
$$
Then a solution of the Cauchy problem (\ref{classic1}),
(\ref{classic2}) is given by means of the integral
\begin{equation}
\label{integrel} u(t,x) = \int_{0}^{t}v(t, \tau, x)d\tau.
\end{equation}
The formulated statement is known as \textit{Duhamel's principle,}
and the integral in (\ref{integrel}) as \textit{Duhamel's integral.}

A similar statement is valid in the case of the Cauchy problem with
a homogeneous initial condition for a first order inhomogeneous
partial differential equation
$$
\frac{\partial u}{\partial t}(t,x) + Cu(t,x) = f(t,x), \quad t > 0,
\,\,x \in R^{n},
$$
where \,$C=C(x,D_x)$\, is a linear differential operator containing
only spatial derivatives, and with coefficients not depending on
\,$t$ (see \cite{BJS64}).

\section{Generalizations of Duhamel's principle}

In this section we prove abstract fractional generalizations of
Duhamel's principle and discuss some of their applications.

\subsection{Duhamel's principle: $\Lambda=\sum_{k=0}^m\delta_{_{\alpha_k}}$ with $\alpha_{k} = k, \,\,k
= 1, ... , m$.} Suppose the measure $\Lambda$ in \eqref{eq1} has the
form $\Lambda=\sum_{k=0}^m\delta_{_{k}},$ where $\delta_{a}$ denotes
Dirac's delta with mass on $a.$ Suppose also that $f(m,A)=I,$ the
identity operator. Then the Cauchy problem (\ref{eq1}),
(\ref{cauchy1}) takes the form
\begin{equation}
\label{eq2}
 u^{(m)}(t) + \sum_{k=0}^{m-1}f_{k}(A) u^{(k)}(t) =
h(t),\quad t>0,
\end{equation}
\begin{equation}
\label{cauchy2}
 u^{(k)}(0)= \varphi_{k}, \quad k = 0, ..., m-1.
\end{equation}
The operators $f_k(A)=f(k,A), k=0,...,m-1,$ are understood in the
sense of the functional calculus introduced in Section
\ref{fcalculus}. In the following theorem we assume that the
vector-functions $U(t,\tau)$ and $h(t)$  are
$\mathrm{Exp}_{A_{\lambda}, G}(\mathrm{X})$-, or
$\mathrm{Exp}^{'}_{A^{\ast}, G^{\ast}}(\mathrm{X^{\ast}})$-valued.
In this abstract case Duhamel's principle is formulated as follows.

\begin{thm} \label{thinteger}
Let a vector-function \,$U(t,\tau)$\, for all $\tau: 0\le \tau<t $
be a solution of the Cauchy problem for a homogeneous equation
\begin{align}
\label{eq3} &\frac{\partial^{m}U}{\partial t^{m}}(t,\tau) +
\sum_{k=0}^{m-1}f_{k}(A) \frac{\partial^{k}U}{\partial
 t^{k}}(t,\tau) = 0,\quad  t > \tau,
\\
\label{cauchy3} &\frac{\partial^{k}U}{\partial
t^{k}}(t,\tau)|_{t=\tau+0} = 0, \,\,k = 0, ... , m-2,
\\
\label{cauchy31} &\frac{\partial^{m-1}U}{\partial
t^{m-1}}(t,\tau)|_{t=\tau+0} = h(\tau),
\end{align}
where $h(t)$ is a continuous vector-function. Then a solution of the
Cauchy problem for the inhomogeneous equation
\begin{equation}
\label{eq4} u^{(m)}(t) + \sum_{k=0}^{m-1}f_{k}(A) u^{(k)}(t) = h(t),
\end{equation}
\begin{equation}
\label{cauchy4} u^{(k)}(0) = 0, \,\,k = 0, ... , m-1.
\end{equation}
is represented via Duhamel's integral
\begin{equation}
\label{solution1} u(t) = \int_{0}^{t}U(t,\tau)d\tau.
\end{equation}
\end{thm}

{\it Proof.} Obviously \,$u(0) = 0$. Further, for the first order
derivative of $u(t)$, using Lemma \ref{lemma2}, one has
$$
\frac{du}{dt}(t) = U(t,t) + \int_{0}^{t}\frac{\partial U}{\partial
t}(t,\tau)d\tau,
$$
By virtue of (\ref{cauchy3}) the latter implies that
\,$\frac{du}{dt}(0) = 0$. Further, differentiating,
$$
\frac{d^{k}u}{dt^{k}}(t) = \frac{\partial^{k-1} U}{\partial
t^{k-1}}(t,t) + \int_{0}^{t}\frac{\partial^{k} U}{\partial
t^{k}}(t,\tau)d\tau,
$$
which due to condition (\ref{cauchy3}) implies that
$$
\frac{d^{k}u}{dt^{k}}(0) = 0,\,\,k = 2, ... , m-1.
$$
Therefore, the function \,$u(t)$\, in (\ref{solution1}) satisfies
initial conditions (\ref{cauchy4}). Moreover, substituting
(\ref{solution1}) to (\ref{eq4}), and taking into account
\eqref{cauchy31}, we have
\begin{align}
u^{(m)}(t) &+  \sum_{k=0}^{m-1}f_{k}(A)u^{(k)}(t) =
    \frac{d^{m}}{dt^{m}} \int_{0}^{t}U(t,\tau)d\tau +
\sum_{k=0}^{m-1}f_{k}(A)\frac{d^{k}}{dt^{k}}
\int_{0}^{t}U(t,\tau)d\tau \notag\\
&= \frac{\partial^{m-1} U}{\partial t^{m-1}}(t,t) +
\int_{0}^{t}\frac{\partial^{m} U}{\partial t^{m}}(t,\tau)d\tau +
\sum_{k=0}^{m-1}f_{k}(A)\int_{0}^{t}\frac{\partial^{k}U}{\partial
t^{k}}(t,\tau)d\tau \notag \\
&=  h(t) + \int_{0}^{t}\left[\frac{\partial^{m} U}{\partial
t^{m}}(t,\tau) +
\sum_{k=0}^{m-1}f_{k}(A)\frac{\partial^{k}U}{\partial
t^{k}}(t,\tau)\right]d\tau =   h(t).\notag
\end{align}
 Hence, \,$u(t)$\, in
(\ref{solution1}) satisfies equation (\ref{eq4}) as well. \eproof

\begin{rem}
It is not hard to see that Theorem \ref{thinteger} holds with
generic closed operators $B_k$ (with dense domain ${\cal{D}}(B_k)$
and commuting with $\frac{d}{dt}$) instead of $f_k(A).$ In this case
we assume that $h(t) \in X$ and $\frac{\partial^k
U(t,\tau)}{\partial t^k} \in {\cal{D}}(B_k), k=0,...,m-1.$

\end{rem}


\subsection{Fractional Duhamel's principle: $\Lambda=\delta_{\mu}+\lambda$ with $\mu \in(m-1,m]$.}

Let $\Lambda=\delta_{\mu}+\lambda,$ where $\mu$ is a number such
that $m-1 < \mu < m,$ and $\lambda$ is a finite measure with $supp
\, \lambda \subset [0,m-1].$
Consider the operator
\begin{equation} \label{elgen}
\,_{\tau}L^{(\mu, \, \lambda)}[u](t)  \equiv \,_{\tau}D_*^{\mu} u(t)
+
\int_0^{m-1}f(\alpha,A) \,_{\tau}D_*^{\alpha} u(t) \lambda(d\alpha), 
\end{equation}
acting on $m$-times differentiable vector-functions $u(t), \, t \ge
\tau \ge 0.$ If $\tau=0,$ then instead of
$\,_{0}L^{(\mu,\lambda)}[u](t)$ we write $L^{(\mu,\lambda)}[u](t).$
\begin{thm} \label{thm1} Suppose that \,$V(t,\tau),\,\,t \ge \tau \geq 0$, is
a solution of the Cauchy problem for the homogeneous equation
\begin{align}
\label{eq5} \,_{\tau}&L^{(\mu, \, \lambda)}[V(\cdot,  \tau)](t) =0,  
\quad t > \tau,\\
\label{cauchy5}
&\frac{\partial ^{k}V}{\partial
t^{k}}(t,\tau)|_{t=\tau+0} = 0,\,\,\,k = 0, ..., m-2,
\\
\label{cauchy51}
&\frac{\partial ^{m-1}V}{\partial
t^{m-1}}(t,\tau)|_{t=\tau+0} = D_{+}^{m-\mu}h(\tau),
\end{align}
where \,$h(t)$\, is a given vector-function. Then Duhamel's integral
\begin{equation}
\label{solution2} u(t) = \int_{0}^{t}V(t,\tau)d\tau
\end{equation}
solves the Cauchy problem for the inhomogeneous equation
\begin{equation}
\label{eq51} L^{(\mu, \, \lambda)}[u](t) = h(t),\,\,\,t > 0,
\end{equation}
 with the homogeneous Cauchy conditions
\begin{equation}
\label{cauchy52} u^{(k)}(0)= 0, \quad k = 0, ..., m-1.
\end{equation}
\end{thm}

{\it Proof.} First notice that since  \,$m-1<\mu < m,$ and therefore
\,$0<m-\mu <1,$\, due to Lemma \ref{lemma1}, the equation
$J^{m-\mu}g(t)=h(t)$ has a unique solution
\begin{equation}
\label{sol6} g(t)= D_{+}^{m - \mu}h(t).
\end{equation}
Let $V(t,\tau)$ as a function of the variable $t$ be a solution to
Cauchy problem (\ref{eq5})-(\ref{cauchy51}) for any fixed $\tau.$ We
verify that \,$u(t)=\int_{0}^t V(t,\tau)d\tau$\, satisfies equation
(\ref{eq51}), and conditions \eqref{cauchy52}. Splitting the
interval $(0,m-1]$ into subintervals $(0,1], ..., (m-2,m-1],$ we
have
\begin{align}
&L^{(\mu, \, \lambda)}[u](t) = D_{\ast}^{\mu}u(t) +
\sum_{k=1}^{m-1}\int_{k-1}^{k}f(\alpha,A)D_{\ast}^{\alpha}u(t)
\lambda(d \alpha).
\label{1ut} 
\end{align}

If $\alpha \in (k-1,k],$ where $k=1,...,m-1$ using the definition of
$D_{\ast}^{\alpha},$
 we have
\begin{align}
D_{\ast}^{\alpha}u(t)&=
\frac{1}{\Gamma(k-\alpha)}\int_0^t(t-s)^{k-\alpha-1}\frac{d^k}{ds^k}\int_0^s
V(s,\tau)d\tau ds. \label{lambda0}
\end{align}
Lemma \ref{lemma2} and conditions (\ref{cauchy5}) 
imply that
\begin{equation}
\label{relation5} \frac{d^{k}}{ds^{k}}\int_{0}^{s} V(s,\tau)d\tau =
\int_{0}^{s}\frac{\partial^{k}}{\partial s^{k}} V(s,\tau)d\tau, \, k
=1,...,m-1.
\end{equation}
Hence, for $\alpha \in (k-1,k], \, k=1,...,m-1,$
\begin{align}
D_{\ast}^{\alpha}u(t)&= \int_0^t
\frac{1}{\Gamma(k-\alpha)}\int_{\tau}^t
(t-s)^{k-\alpha-1}\frac{\partial^k}{\partial s^k}V(s,\tau)ds d\tau.
\label{lambda0}
\end{align}
Again due to Lemma \ref{lemma2} and condition 
(\ref{cauchy51}),
\begin{align}
\frac{d^{m}}{ds^{m}}\int_{0}^{s} V(s,\tau)d\tau &= \frac{\partial
^{m-1}}{\partial s^{m-1}}V(s,\tau)_{_{|\tau=s}} +
\int_{0}^{s}\frac{\partial^{m}}{\partial s^{m}} V(s,\tau)d\tau
\notag
\\ &= D_{+}^{m-\mu}h(s) +
 \int_{0}^{s}\frac{\partial^{m}}{\partial
s^{m}} V(s,\tau)d\tau. \label{relation6}
\end{align}
Therefore the first term on the right hand side of \eqref{1ut} takes
the form
\begin{align}
D_{\ast}^{\mu}u(t)&=
\frac{1}{\Gamma(m-\mu)}\int_{0}^{t}(t-s)^{m-\mu-1}
\frac{d^{m}}{ds^{m}}\int_{0}^{s} V(s,\tau)d\tau ds \notag
\\
&= \frac{1}{\Gamma(m-\mu)}\int_{0}^{t}(t-s)^{m-\mu-1}
\Big(D_{+}^{m-\mu}h(s) +
 \int_{0}^{s}\frac{\partial^{m}}{\partial
s^{m}} V(s,\tau)d\tau\Big)ds. \label{delta0}
\end{align}
Further, put $g(t)=D_{+}^{m-\mu}h(t).$ Then by virtue of
\eqref{sol6},
\begin{equation} \label{ht}
\frac{1}{\Gamma(m-\mu)} \int_{0}^{t} (t-s)^{m-\mu-1}
D_{+}^{m-\mu}h(s)= J^{m-\mu} g(t)= h(t).
\end{equation}
Now equations (\ref{1ut}), (\ref{lambda0}), \eqref{delta0}, and
(\ref{ht}) imply that
\begin{align}
L^{(\mu, \, \lambda)}&[u](t) = h(t) +
\frac{1}{\Gamma(m-\mu)}\int_{0}^{t}(t-s)^{m-\mu-1} \int_{0}^{s}
\frac{\partial^{m}}{\partial
s^{m}}V(s,\tau)d\tau\,ds \notag \\
&+ \sum_{k=1}^{m-1} \int_{k-1}^k f(\alpha,A)
\frac{1}{\Gamma(k-\alpha)}\int_{0}^{t}(t-s)^{k-\alpha-1}
\int_{0}^{s}\frac{\partial^{k}}{\partial s^{k}}
V(s,\tau)d\tau\, ds \lambda(d\alpha). 
\label{1ut2} 
\end{align}
Changing the order of integration (Fubini is allowed) in
(\ref{1ut2}) we get
\begin{align*}
&L^{(\mu, \, \lambda)}[u](t) = h(t) +
\int_{0}^{t}\int_{\tau}^{t}\frac{1}{\Gamma(m-\mu)}(t-s)^{m-\mu-1}
 \frac{\partial^{m}}{\partial
s^{m}}V(s,\tau)ds\,d\tau \\
&+ \sum_{k=1}^{m-1} \int_0^t \int_{k-1}^k f(\alpha,A)
\int_{\tau}^{t}\frac{1}{\Gamma(k-\alpha)}(t-s)^{k-\alpha-1}
\frac{\partial^{k}}{\partial s^{k}} V(s,\tau)ds
\lambda(d\alpha)\,d\tau 
\\
&= h(t) + \int_{0}^{t} \,_{\tau}D_{*}^{\alpha_{m}}V(t,\tau)d\tau\,+
 \int_{0}^{t} \int_0^{m-1} f(\alpha,A) ~ \,_{\tau}D_{*}^{\alpha}
V(t,\tau)\lambda(d\alpha)d\tau\,
\\
& 
=  h(t) + \int_0^t \,_{\tau}L^{(\mu, \, \lambda)}[V(\cdot,\tau)](t)
d\tau = h(t).
\end{align*}
Finally, using the relations (\ref{relation5}) it is not hard to
verify that \,$u(t)$\, in (\ref{solution2}) satisfies initial
conditions (\ref{cauchy52}) as well. \eproof

If the vector-function \,$h$\, satisfies the additional condition
\,$h(0) = 0$\, then condition (\ref{cauchy51}) in accordance with
the relationship (\ref{relation2}) can be replaced by
$$
\frac{\partial ^{m-1}V}{\partial t^{m-1}}(t,\tau)|_{t=\tau} =
D_{*}^{m-\mu} h(\tau),
$$
with the Caputo-Djrbashian derivative \,$D_{*}^{m-\mu}$\, of order
\,$m - \mu$. As a consequence the formulation of the fractional
Duhamel's principle takes the form:

\begin{thm} \label{thm2} Suppose that for all $\tau: ~ 0<\tau<t$ a function \,$V(t,\tau)$, is
a solution to the Cauchy problem for the homogeneous equation
\begin{align*}
\,_{\tau}&L^{(\mu, \, \lambda)}[V(\cdot,\tau)](t) = 0,\quad t > \tau, \\
& \frac{\partial ^{k}V}{\partial t^{k}}(t,\tau)|_{t=\tau+0} =
0,\,\,\,k = 0, ... , m-2, \\
& \frac{\partial ^{m-1}V}{\partial t^{m-1}}(t,\tau)|_{t=\tau+0} =
D_{*}^{m-\mu}h(\tau),
\end{align*}
where \,$h(t)$\, is a given vector-function such that \,$h(0) = 0$.
Then
$$
u(t) = \int_{0}^{t}V(t,\tau)d\tau
$$
is a solution of the Cauchy problem for the inhomogeneous equation
$$
L^{(\mu, \, \lambda)}[u](t) = h(t), \quad t > 0,
$$
 with the homogeneous Cauchy conditions
$$
u^{(k)}(0)= 0, \quad k = 0, ..., m-1.
$$
\end{thm}

\begin{rem}
\begin{enumerate}
\item
 Lemma \ref{lemma1} can be extended to absolutely continuous functions $h(t)$ with
 an appropriate meaning of solution in equation \eqref{Abel} (see, \cite{SKM93}) .
 It is also known \cite{SKM93} that the fractional derivative
$D_{+}^{k-\mu}h(t),
  \,\,k-1 < \mu < k,\,\,k = 1, ... , m,$\, exists a.e., if \,$h(t)$
   is an absolutely continuous function on \,$[0;T]$ for any $T>0$.
   These two facts imply that generalized Duhamel's principles
   proved above hold true for absolutely continuous functions $h(t).$
\item
     In Theorems \ref{thm1} and \ref{thm2} we assumed that
     $f(\mu,A)$ is the identity operator (see equation
     \eqref{elgen}). In the general case, with
     appropriate selection of $G$ we can assume that the inverse operator
     $[f(\mu,A)]^{-1}$ exists. Then with the condition
     $$\frac{\partial^{m-1}V(t,\tau)}{\partial t^{m-1}}\big|_{{t=\tau+}}=[f(\mu,A)]^{-1}D_+^{m-\mu}h(\tau)
     $$ instead of \eqref{cauchy51}, Theorems \ref{thm1} and
     \ref{thm2} remain valid.
\end{enumerate}
\end{rem}

\subsection{Fractional Duhamels's principle with Riemann-Liuoville derivative}

The operator $\,_{\tau}L^{\Lambda}$ in Theorem \ref{thm1} is defined
via the fractional derivative in the sense of Caputo-Djrbashian. A
fractional generalization of Duhamel's principle is also possible
when this operator is defined via the Riemann-Liouville fractional
derivative. In this section we
briefly discuss this important case proving the corresponding
theorem in the simple particular case
$\,_{\tau}L[u](t)=\,_{\tau}D^{\alpha}_+u(t)+Bu(t),$ where
$0<\alpha<1,$ and $B$ is a closed operator with a domain
${\cal{D}}(B)$ dense in $X.$ The general case can be treated in a
similar manner.

\begin{thm} \label{thm4} Suppose that \,$V(t,\tau),\,\,t \ge \tau \geq 0$, is
a solution of the Cauchy type problem for the homogeneous equation
\begin{align}
\label{eq5rl} \,_{\tau}D^{\alpha}_+V(t,\tau)+BV(t,\tau) =0,  
\quad t > \tau,
\end{align}
\begin{equation}
\label{cauchy5rl} \,_{\tau}J^{1-\alpha}V(t,\tau)_{|_{t=\tau+}}
=h(\tau),
\end{equation}
where \,$h(\tau), \, \tau >0,$\, is a continuous vector-function.
Then Duhamel's integral
\begin{equation} \label{solution2rl}
u(t) = \int_{0}^{t}V(t,\tau)d\tau
\end{equation}
solves the Cauchy type problem for the inhomogeneous equation
\begin{equation}
\label{eq51rl} D_+^{\alpha}u(t)+Bu(t) = h(t),\,\,\,t > 0,
\end{equation}
 with the homogeneous initial condition
$J^{1-\alpha}u(0)= 0.$ 
\end{thm}

\textit{Proof.} Let $V(t,\tau)$ satisfy the conditions of the
theorem. Then for Duhamel's integral \eqref{solution2rl}, by virtue
of Lemma \ref{lemma2}, we have
\begin{align}
D_+^{\alpha}u(t)&+Bu(t)=
\frac{1}{\Gamma(1-\alpha)}\frac{d}{dx}\int_0^t
\frac{\int_0^sV(s,\tau)ds}{(t-s)^{\alpha}}d\tau +\int_0^tBV(t,\tau)
d\tau  \notag \\
&=\frac{d}{dx}\int_0^t 
\,_{\tau}J^{1-\alpha}V(t,\tau)d\tau +\int_0^tBV(t,\tau)
d\tau  \notag \\
&=\,_{\tau}J^{1-\alpha}V(t,\tau)_{|_{\tau=t}} +\int_0^t
[\,_{\tau}D^{\alpha}_+V(t,\tau)+BV(t,\tau)]d\tau = h(t).
\end{align}
On the other hand, changing the order of integration and using the
mean value theorem, we obtain
\begin{align} \label{h0}
\|J^{1-\alpha} u(t)\|=\|\int_0^t
\,_{\tau}J^{1-\alpha}V(t,\tau)d\tau\| \le t
\|\,_{\tau}J^{1-\alpha}V(t,\tau)\|.
\end{align}
Condition \eqref{cauchy5rl} implies that $\lim_{t\to
0+}\,_{\tau}J^{1-\alpha}V(t,\tau)=h(0)$ in the norm of $X.$ It
follows from \eqref{h0} that $\lim_{t\to 0+}J^{1-\alpha} u(t)=0$ in
the norm of $X.$ \eproof

\subsection{Applications. Existence and uniqueness theorems}

Theorems 3.1 and 3.3 lead to generalization of the existence and
uniqueness results obtained in papers \cite{GLU00,Sa06} for the
abstract Cauchy problems. 
Let $L^{\Lambda}$ be the distributed fractional order abstract
differential-operator 
defined in \eqref{elgen} with $\tau=0,$ and via the characteristic
function
$$
\Delta(s, z) = s^{\mu} + \int_0^{m-1}f(\alpha,z)s^{\alpha} d \lambda, 
$$
where $\mu \in (m-1,m],$ $\lambda$ is a finite measure with $supp \,
\lambda \subset [0,m-1],$ and $f(\alpha,z)$ is a function continuous
in $\alpha$ and analytic in $z \in G \subset \cal{C}.$ Denote by
$\hat{v}(s)={\cal{L}}[v](s)$ the Laplace transform of a
vector-function $v(t),$ namely
\[
{\cal{L}}[v](s) = \int_0^{\infty} e^{-st}v(t)dt, ~ s>s_0,
\]
where $s_0 \ge 0$ is a real number. It is not hard to verify that if
$v(t) \in Exp_{A,G}(X)$ for each $t \ge 0$ and satisfies the
condition $\|v(t)\| \le C e^{\gamma t}, t\ge 0,$ with some constants
$C>0$ and $\gamma$, then $\hat{v}(s)$ exists and $$\|A^k
\hat{v}(s)\| \le \frac{C_s}{s-\gamma} \nu^k, ~ ~ s>\gamma,$$
implying $\hat{v}(s) \in Exp_{A,G}(X)$ for each fixed $s > \gamma.$
The lemma below gives a formal representation formula for a
solution of the general abstract Cauchy problem 
\begin{align}
&L^{\Lambda}[u](t)=h(t), \, t>0, \label{gen}\\
&u^{(k)}(0+)=\varphi_k, \, k=0,...,m-1. \label{gencon}
\end{align}
\begin{lem}
\label{repr} Let $ c_{\beta}(t,z) =
{\cal{L}}^{-1}[\frac{s^{\beta}}{\Delta(s, z)}](t), ~ z \in G \subset
{\cal{C}},$ where ${\cal{L}}^{-1}$ stands for the inverse Laplace
transform, and
\[
S_{k}(t,z) = c_{\mu - k - 1}(t,z) + 
\int_{k}^{m-1}f(\alpha,z)c_{\alpha - k - 1}(t,z)\lambda(d\alpha).
\]
Then $u_k(t)=S_k(t,A)\varphi_k$ solves the Cauchy problem
\[
L[u]=0, ~~~ u^{j}(0)=\delta_{j,k}\varphi_j, ~ j=0,...,m-1,
\]
where $\delta_{j,k}=1$ if $j=k,$ and $\delta_{j,k}=0,$ if $j \neq
k.$
\end{lem}

\begin{cor}
\label{cor} Let $S_k(t,A), k=0,...,m-1,$ be the collection of
solution operators with the symbols $S_k(t,z)$ defined in Lemma
\ref{repr}. Then the solution of the Cauchy problem
\begin{equation}\label{coreq1}
 L[u]=0, ~~~ u^{j}(0)=\varphi_j, ~ j=0,...,m-1.
\end{equation}
is given by the following representation formula
\begin{equation}
\label{formal} u(t) = \sum_{k=0}^{m-1} S_k(t,A) \varphi_k.
\end{equation}
\end{cor}

{\em Proof of Lamma.} Applying formula (\ref{laplace}) we have
\[
{\cal{L}}[L [u]](s) = s^{\mu} \hat{u}(s) - \sum_{i=0}^{m-1} u^{i}(0)
s^{\mu -i-1} + \]
\[\sum_{k=1}^{m-1} \int_{k-1}^{k} f(\alpha, A) (s^{\alpha}
\hat{u}(s) - \sum_{j=0}^{k-1} u^j(0) s^{\alpha -j -1})
\lambda(d\alpha)=0. 
\]
Due to the initial conditions $u^{j}(0)=\delta_{j,k}\varphi_j, ~
j=0,...,m-1,$ the latter reduces to
\[
\Delta(s, z) \hat{u}(s)= \varphi_{k} \left(s^{\mu-k-1} +
\int_{k}^{m-1} f(\alpha,z) s^{\alpha-k-1}\lambda(d\alpha)\right).
\]
Now it is easy to see that the solution in this case is represented
$u_k=S_k(t,A)\varphi_k.$ \eproof

\begin{rem} \label{rem3}
\begin{enumerate}
\item Corollary \ref{cor} can easily be extended to the operator
$\,_{\tau}L$ in \eqref{coreq1} as well with the initial conditions
$u^j(\tau)=\varphi_j,$ maintaining the shift $t^{'}=t-\tau;$
\item
A particular case of Lemma \ref{repr} when
$\Lambda=\sum_{k=0}^m\delta_{_{\alpha_k}}, \, k-1<\alpha_k <k$, is
proved in \cite{GLU00}.
\end{enumerate}
\end{rem}

Further, denote by \,$C^{(m)}[t > 0;\,Exp_{A,\,G}(X)]$ \, and by
\,$AC[t > 0;\,Exp_{A,\,G}(X)]$ \, the space of \, $m - $ times
continuously differentiable functions and the space of absolutely
continuous functions on \,$(0; +\infty)$ \, with values ranging in
the space \,$Exp_{A,\,G}(X),$\, respectively. A vector-function
\,$u(t) \in C^{(m)}[t>0;\,Exp_{A,\,G}(X)]$\,
 $\cap \,\, C^{(m-1)}[t \geq 0;\, Exp_{A,\,G}(X)]$\, is called a
solution of the problem  (\ref{gen}), (\ref{gencon})  if it
satisfies the equation (\ref{gen}) and the initial conditions
(\ref{gencon}) in the topology of $Exp_{A,\,G}(X).$

Theorem \ref{thm1} and Corollary \ref{cor} imply the following
results.

\begin{thm} \label{thm11}
Let \,$\varphi_{k} \in Exp_{A,\,G}(X)$,\,\,$k = 0,\,. . .\,,\, m-1$,
\quad $h(t) \in AC[0 \le t \le T;\, Exp_{A,\,G}(X)]$ for any $T>0,$
and \,$D_{+}^{m - \alpha_{m}}h(t) \in C[0 \le t \le T;\,
Exp_{A,\,G}(X)]$. Then the Cauchy problem (\ref{gen}),
(\ref{gencon}) has a unique solution. This solution is given by 
\begin{equation} \label{repr1}
u(t) = \sum_{k=0}^{m-1}S_{k}(t,A)\varphi_{k} +
\int_{0}^{t}S_{m-1}(t-\tau,A)D_{+}^{m - \mu} h(\tau)d\tau.
\end{equation}
\end{thm}

\textit{Proof.} We split the Cauchy problem
\eqref{gen},\eqref{gencon} into two Cauchy problems
\begin{align}
&L^{\Lambda}[u](t)=0, \, t>0, \label{gen1}\\
&u^{(k)}(0+)=\varphi_k, \, k=0,...,m-1, \label{gencon1}
\end{align}
and
\begin{align}
&L^{\Lambda}[v](t)=h(t), \, t>0, \label{gen2}\\
&v^{(k)}(0+)=0, \, k=0,...,m-1. \label{gencon2}
\end{align}
Due to Corollary \ref{cor} the unique solution to
\eqref{gen1},\eqref{gencon1} is given by
\begin{equation} \label{repr11}
u(t) = \sum_{k=0}^{m-1}S_{k}(t,A)\varphi_{k}.
\end{equation}
For the Cauchy problem \eqref{gen2},\eqref{gencon2}, in accordance
with the fractional Duhamel's principle (Theorem \ref{thm1}), it
suffices to solve the Cauchy problem for the homogeneous equation:
\begin{align}
\,_{\tau}&L^{\Lambda}[V(t,\tau)](t)=0, \, t>\tau, \label{gen3}\\
&\frac{\partial^k V(t,\tau)}{\partial t^k}_{\big|{t=\tau+}}=0, \,
k=0,...,m-2, \label{gencon3}\\
&\frac{\partial^{m-1} V(t,\tau)}{\partial
t^{m-1}}_{\big|t=\tau+}=D_+^{m-\mu}h(\tau). \label{gencon33}
\end{align}
The solution of this problem, again using Corollary \ref{cor} (with
the note in Remark \ref{rem3}), has the representation
\begin{equation}
V(t,\tau)=S_{m-1}(t-\tau,A)D_{+}^{m - \mu} h(\tau).
\end{equation}
Thus, Duhamel's integral of the latter and representation
\eqref{repr11} lead to formula \eqref{repr1}. The uniqueness of a
solution also follows from the obtained representation \eqref{repr1}
(see \cite{Um98}). \eproof

The duality immediately implies the following theorem.
\begin{thm} \label{thm12}
Let \,$\varphi^{\ast}_{k} \in
Exp^{'}_{A^{\ast},\,G^{\ast}}(X^{\ast})$,\,\,$k = 0,\,. . .\,,\,
m-1$, \quad $h^{\ast}(t) \in AC[t \le t \le T;\,
Exp^{'}_{A^{\ast},\,G^{\ast}}(X^{\ast})]$\, and \,$D_{+}^{m -
\alpha_{m}}h^{\ast}(t) \in C[t \le t \leq T;\,
Exp^{'}_{A^{\ast},\,G^{\ast}}(X^{\ast})]$. Assume also that
$Exp_{A,\,G}(X)$ is dense in $X.$ Then the Cauchy problem
(\ref{gen}), (\ref{gencon}) (with $A$ switched to $A^{\ast}$) is
meaningful and has a unique weak solution. This solution is given by
$$
u^{\ast}(t) = \sum_{k=0}^{m-1}S_{k}(t,A^{\ast})\varphi^{\ast}_{k} +
\int_{0}^{t}S_{m-1}(t-\tau,A^{\ast})D_{+}^{m - \mu}
h^{\ast}(\tau)d\tau.
$$
\end{thm}

Assume that $Exp_{A,\,G}(X)$ is densely embedded into $X.$ Besides,
let the solution operators $S_k(t,A)$ for each $k=0,...,m-1,$
satisfy the estimates
\begin{equation} \label{EST}
\|S_k(t,A) \varphi\| \le C \|\varphi\|, ~~ \forall ~ t \in [0,T],
\end{equation}
where $\varphi \in Exp_{A,\,G}(X),$ and $C>0$ does not depend on
$\varphi.$ Then there exists a unique closure $\bar{S}_k (t)$ to $X$
of the operator $S_k(t,A)$ which satisfies the estimate $\|\bar{S}_k
(t) u\| \le\|u\| $ for all $u \in X.$ Using the standard technique
of closure (see \cite{Um97,Um98}), we can prove the following
theorem.

\begin{thm} \label{thm13}
Let \,$\varphi_{k} \in X$,\,\,$k = 0,\,. . .\,,\, m-1$, \quad $h(t)
\in AC[0 \le t \le T;\, X]$ for any $T>0,$ and \,$D_{+}^{m -
\alpha_{m}}h(t) \in C[0 \le t \le T;\, X]$. Further let
$Exp_{A,\,G}(X)$ be densely embedded into $X,$ and the estimates
(\ref{EST}) hold for solution operators $S_k(t,A), k=0,...,m-1.$
Then the Cauchy problem (\ref{gen}), (\ref{gencon}) has a unique
solution $u(t) \in C^{m}[0 < t \le T;\, X]$. This solution is given
by 
$$
u(t) = \sum_{k=0}^{m-1}\bar{S}_{k} (t) \varphi_{k} + \int_{0}^{t}
\bar{S}_{m-1}(t-\tau)D_{+}^{m - \mu} h(\tau)d\tau.
$$
\end{thm}

\end{document}